\documentclass[runningheads]{llncs}
\usepackage{graphicx}
\usepackage{amsmath,amssymb,euscript,yfonts,psfrag,latexsym,dsfont}
\newcommand{\mE}{{\mathbb E}}

\newcommand{\cP}{{\mathcal P}}
\newcommand{\cQ}{{\mathcal Q}}

\newcommand{\argmin}{\operatorname{argmin}}
\newcommand{\Proj}{\operatorname{Proj}}

\begin{document}
\title{Multi-marginal Schr\"odinger bridges\thanks{Partial support was provided by NSF under grants 1665031, 1807664, 1839441 and 1901599, and by AFOSR under grant FA9550-17-1-0435.}}
%
%
\author{Yongxin Chen\inst{1} \and
Giovanni Conforti\inst{2} \and
Tryphon T. Georgiou \inst{3}
\and Luigia Ripani\inst{3}}
\authorrunning{Chen et al.}
%
\institute{Georgia Institute of Technology, Atlanta, GA 30332, USA \and
\'Ecole Polytechnique, Route de Saclay, 91128, Palaiseau Cedex, France
 \and
University of California, Irvine, CA 92697, USA}
\maketitle              
\begin{abstract} We consider the problem to identify the most likely flow in phase space, of (inertial) particles under stochastic forcing, that is in agreement with spatial (marginal) distributions that are specified at a set of points in time. The question raised generalizes the classical Schr\"odinger Bridge Problem (SBP) which seeks to interpolate two specified end-point marginal distributions of overdamped particles driven by stochastic excitation. While we restrict our analysis to second-order dynamics for the particles, the data represents {\em partial} (i.e., only positional) information on the flow at {\em multiple} time-points. The solution sought, as in SBP, represents a probability law on the space of paths this closest to a uniform prior while consistent with the given marginals.
We approach this problem as an optimal control problem to minimize an action integral {\em a la} Benamou-Brenier, and derive a 
 time-symmetric formulation that includes a Fisher information term on the velocity field. We underscore the relation of our problem to recent measure-valued splines in Wasserstein space, which is akin to that between SBP and Optimal Mass Transport (OMT). The connection between the two provides a Sinkhorn-like approach to computing measure-valued splines.
We envision that interpolation between measures as sought herein will have a wide range of applications in signal/images processing as well as in data science in cases where data have a temporal dimension. 

\end{abstract}
\section{Introduction}
In 1931/32, in an attempt to gain insights into the stochastic nature of quantum mechanics, Schr\"odinger \cite{Sch31,Sch32} raised the following question regarding a system of a large number of classical independent identically distributed (i.i.d.) Brownian particles. He hypothesized that this ``cloud'' of particles is observed to have (empirical) distributions $\rho_0(x_0)$ and $\rho_1(x_1)$ at two points in time $t_0=0$ and $t_1=1$, respectively, and
further that $\rho_1(x_1)$ differs from what is dictated by the law of large numbers, i.e., that
    \[
       \rho_1(x_1)\neq \int_{{\mathbb R}^n} q(t_0,x_0,t_1,x_1)\rho_0(x_0)dx_0,
    \]
where
    \[
        q(s,x,t,y)=(2\pi)^{-n/2}(t-s)^{-n/2}\exp\left(-\frac12\frac{\|x-y\|^2}{(t-s)}\right)
    \]
denotes the Brownian transition probability kernel. Schr\"odinger then sought to find the ``most likely'' evolution for the cloud of particles to have transitioned from $\rho_0$ to $\rho_1$. 
In the language or large deviation theory (which was not in place at the time),   Schr\"odinger's question amounts to seeking a probability law on the path space that is in agreement with the two marginals while being the closest to the Brownian prior in the sense of relative entropy \cite{Fol88}. The solution is known as the Schr\"odinger bridge since the law ``bridges'' the two given end-point marginal distributions.

Renewed interest in the Schr\"odinger Bridge Problem (SBP) has been fueled by its connections to the Monge-Kantorovic Optimal Mass Transport (OMT) and a wide range aplications in image analysis, stochastic control, and physics \cite{Leo12,Leo13,CheGeoPav14e,GenLeoRip15,CheGeoPav14a,conforti2018second,peyre2018computational}. More specifically, SBP, seen as a suitable regularization of OMT, provides a natural model for uncertainty in the transport of distributions as well as a valuable computational tool for interpolating distributional data.

In this work we consider a natural generalization of the Schr\"odinger bridge theory to address the situation where the data consist of possible partial marginal distributions at various points in time.
Thus, we postulate a similar experiment with stochastic particles. However, in contrast to the standard SB theory, we concieve these particles to obey second order stochastic differential equations with Brownian stochastic forcing that accounts for random acceleration along trajectories in phase space (see Section \ref{sec:SB}). This new setting connects with recent results on measure-valued splines \cite{CheConGeo18,BenGalVia18,CheKar18}, which are general notions of splines on the Wasserstein space of measures. In this short paper, we provide a summary of the theory. A more detailed account will appear in a forthcoming publication (in preparation).

\section{Multi-marginal Schr\"odinger bridges for inertial particles}\label{sec:SB}
Suppose we are given a large number of independent inertial particles driven by white noise, that is, they follow the dynamics
	\begin{subequations}\label{eq:sys}
	\begin{eqnarray}
	dx &=& v dt,\\
	dv &=& dw,
	\end{eqnarray}
	\end{subequations}
where $dw$ denotes standard Brownian motion and $(x,\,v)$ is the phase space ($x$ denoting position and $v$ velocity).
The flow of probability densities $\mu_t(x,v)$ in phase space obeys the Fokker-Planck equation
	\begin{equation}
	\frac{\partial \mu}{\partial t} +v\cdot\nabla_x\mu -\frac{1}{2}\Delta_v \mu = 0,
	\end{equation}
with initial condition $\mu_0$ at time $t=0$. The law of large numbers dictates that, when the number of particles is large enough, the distributions of the particles will be closed to $\mu_t$. The data of the problem we consider will, as in the standard SB problem, be inconsistent with the Fokker-Planck equation, and viewed as a ``atypical/rare'' event.
In standard SB setting, where only two end-point marginals  are specified, the ``most likely'' evolution amounts to an adjustment of the Fokker-Planck equation by adding a suitable drift term to match the two marginals.
For inertial particles driven by white noise, the generator is hypoelliptic \cite{Hor67} and the SB theory carries over to matching marginals in phase space. 

Throughout, we suppose that we have only access to position $x$ and that empirical marginals $\rho_0 = \Proj_x \mu_0, \rho_1 = \Proj_x \mu_1$ represent projections, accordingly. The extra degree of freedom, since $\mu$'s are partially specified, make the corresponing multi-marginal problem nontrivial. Specifically, we seek the most like paths the particles have taken that match positional distributions $\rho_0, \rho_1, \cdots, \rho_N$ at times $0= t_0<t_1<\cdots<t_N=1$. To this end, we let $\cQ$ be the law of \eqref{eq:sys} (on path space) and we let $\cP$ be the any other law. We seek the minimizer of
	\begin{equation}
	H(\cP,\,\cQ) = \int d\cP \log\frac{d\cP}{d\cQ}
	\end{equation}
over all laws $\cP$ that are consistent with the marginals $\rho_0, \rho_1, \cdots, \rho_N$.

To guarantee the boundedness of the relative entropy $H$ between $\cP$ and the prior process $\cQ$, $\cP$ has to be of the form
	\begin{subequations}\label{eq:drift}
	\begin{eqnarray}
	dx &=& v dt,\\
	dv &=& a dt + dw,
	\end{eqnarray}
	\end{subequations}
where $a$ is a suitable drift that may depend on the current and past values of the process state. 
Invoking Girsanov's theorem \cite{Gir60,IkeWat14,Dai91,CheGeoPav15b}, we obtain that
	\[
		H(\cP,\,\cQ) = \mE_{\cP} \left\{\int_0^1\|a(t)\|^2 dt\right\}=\int_0^1 \int \|a\|^2 \mu dx dv dt.
	\]
Thus, we arrive at the optimal control formulation
	\begin{subequations}\label{eq:control}
	\begin{eqnarray}
	\min&&\mE \left\{\int_0^1\|a(t)\|^2 dt\right\}, \label{eq:control1}
	\\&&
	dx = v dt, \label{eq:control2}
	\\&&
	dv = a dt + dw, \label{eq:control3}
	\\&& x(t_i)~\sim~ \rho_i,~~i=0, 1, \ldots, N. \label{eq:control4}
	\end{eqnarray}
	\end{subequations}
The difference to the standard Schr\"odinger bridge problem lies in the constraint \eqref{eq:control4}; {\em multiple} marginals are {\em partially} specified. The existence and uniqueness of the solution follow from the fact that it is a strongly convex optimization problem on path space measures. The argument is similar to the standard argument in SB theory \cite{Fol88,Leo12,Leo13,CheGeoPav15a} and will be presented in an extended version of the paper. By utilizing the Fokker-Planck equation, we can rewrite the above as
	\begin{subequations}\label{eq:SBSO}
	\begin{eqnarray}
	\min&&\int_0^1 \int \|a\|^2 \mu dx dv dt,\label{eq:SBSO1}
	\\&& \frac{\partial \mu}{\partial t} +v\cdot\nabla_x\mu+\nabla_v\cdot(a \mu) -\frac{1}{2}\Delta_v \mu = 0,\label{eq:SBSO2}
	\\&& \int \mu_{t_i}(x,v)dv = \rho_i(x),~~i=0, 1, \ldots, N.\label{eq:SBSO3}
	\end{eqnarray}
	\end{subequations}
Here $\int \mu_{t_i}(x,v)dv = \rho_i(x)$ since $\Proj_x(\mu_{t_i}) = \rho_i$.

Now let $\hat a = a - \frac{1}{2}\nabla_v \log \mu$, then the diffusion term in \eqref{eq:SBSO2} can be absorbed into the convection terms, and then the cost becomes
	\begin{eqnarray*}
	\int_0^1 \int \|a\|^2 \mu dx dv dt &&=
	\int_0^1 \int \|\hat{a}+\frac{1}{2}\nabla_v \log \mu\|^2 \mu dx dv dt
	\\&&=\int_0^1 \int \left\{\|\hat{a}\|^2 \mu +\frac{1}{4}\|\nabla_v \log \mu\|^2 \mu\right\}dx dv dt
	\\&&+ \int_0^1 \int \langle \hat a, \nabla_v \log \mu\rangle\mu dx dv dt.
	\end{eqnarray*}
Direct calculation yields 
	\begin{equation}
	\int_0^1 \int \langle \hat a, \nabla_v \log \mu\rangle\mu dx dv dt
	= \int \{\mu_1 \log \mu_1-\mu_0\log\mu_0\} dx dv,
	\end{equation}
which only depends on the two end distributions.
Thus, we need to consider 
	\begin{subequations}\label{eq:SBSOS}
	\begin{eqnarray}
	\hspace{-0.45cm}\min\!&&\int_0^1 \hspace*{-5pt}\int \left\{\|\hat{a}\|^2 \mu \!+\!\frac{1}{4}\|\nabla_v \log \mu\|^2 \mu\right\}dx dv dt
	\!+\!\int \{\mu_1 \log \mu_1\!-\!\mu_0\log\mu_0\} dx dv,
	\\&& \frac{\partial \mu}{\partial t} +v\cdot\nabla_x\mu+\nabla_v\cdot(\hat{a} \mu) = 0,
	\\&& \int \mu_{t_i}(x,v)dv = \rho_i(x),~~i=0, 1, \ldots, N.
	\end{eqnarray}
	\end{subequations}
A similar formulation with a Fisher information term has been studied in \cite{Yas81,CheGeoPav14e,GenLeoRip15,LiYinOsh18} for standard Schr\"odinger bridge problems. However, in the standard setting, the term $\int \{\mu_1 \log \mu_1\!-\!\mu_0\log\mu_0\} dx dv$ can be dropped since the full-state marginal distributions are already specified. It important to note that, compared to \eqref{eq:SBSO}, \eqref{eq:SBSOS} is time symmetric. 


\section{Connections to measure-valued splines}\label{sec:spline}
A variational formulation \cite{Hol57} of splines going through $\{x_1, x_2, \ldots, x_N\}$ in Euclidean space is given by
	\begin{eqnarray*}
		\min &&\int_0^1 \|\ddot x\|^2 dt
		\\&&
		x(t_i) = x_i, ~~i=0, 1, \ldots, N,
	\end{eqnarray*}
where the minimization is taken over all twice-differentiable trajectories that satisfy the constraints. This formation has been generalized to the Wasserstein space of measures \cite{CheConGeo18,BenGalVia18}. In particular, the fluid dynamic formulation for measure-valued splines in \cite{CheConGeo18} with marginals $\rho_0, \rho_1, \cdots, \rho_N$ at $0= t_0<t_1<\cdots<t_N=1$ reads
	\begin{subequations}
	\begin{eqnarray}
	\inf&&\int_0^1 \int \|a\|^2 \mu dx dv dt,
	\\&& \frac{\partial \mu}{\partial t} +v\cdot\nabla_x\mu+\nabla_v\cdot(a \mu) = 0,
	\\&& \int \mu_{t_i}(x,v)dv = \rho_i(x),~~i=0, 1, \ldots, N.
	\end{eqnarray}
	\end{subequations}
We note that the above formulation is almost the same as \eqref{eq:SBSO} except for a missing diffusion term in the constraint. This resembles the relation between standard Schr\"odinger bridges and optimal mass transport \cite{Leo12,Leo13,CheGeoPav15a}.

Indeed, a zero-noise limit argument follows. If we replace the dynamics \eqref{eq:sys} of the inertial particles by
	\begin{subequations}
	\begin{eqnarray}
	dx &=& v dt,\\
	dv &=& \sqrt{\epsilon}dw,
	\end{eqnarray}
	\end{subequations}
then the multi-marginal SB problem becomes
	\begin{subequations}\label{eq:SBG}
	\begin{eqnarray}
	\min &&\int_0^1 \int \|a\|^2 \mu dx dv dt,
	\\&& \frac{\partial \mu}{\partial t} +v\cdot\nabla_x\mu+\nabla_v\cdot(a \mu) -\frac{\epsilon}{2}\Delta_v \mu = 0,
	\\&& \int \mu(t_i, x,v)dv = \rho_i(x),~~i=0, 1, \ldots, N.
	\end{eqnarray}
	\end{subequations}
The ``slowed down'' formulation \eqref{eq:SBG} reduces to \eqref{eq:SBSO} when we take the limit $\epsilon \rightarrow 0$. 	
Therefore, we establish the measure-valued spline as a zero-noise limit of a multi-marginal SB, and SB as a regularized version of measure-valued spline. Rigorous proof of these conclusions will be presented in a forthcoming paper. 
	
\section{Algorithms}
The Sinkhorn \cite{Sin64,Cut13} algorithm is a natural iterative scheme in SB problems. Due to its efficiency and simplicity, it has became a workhorse for data science applications of optimal transport \cite{Cut13}. In this section, we develop a Sinkhorn-type algorithm for multi-marginal Schr\"odinger bridge problems. 
The relation established in Section \ref{sec:spline} implies that the same algorithm can also be used to approximate measure-valued splines. 

Using a measure decomposition argument, the Schr\"odinger problem \eqref{eq:SBG} can be rewritten as 
	\begin{subequations}
	\begin{eqnarray}
	&&\min J(\pi):=\sum_{i=0}^{N-1} KL (\pi_{i,i+1}\mid e^{-C_{i,i+1}/\epsilon})
	\\&&\int \pi_{i,i+1} dx_{i+1}dv_{i+1} = \mu_i,~~i=0,\ldots, N-1
	\\&&\int \pi_{i,i+1} dx_{i}dv_{i} = \mu_{i+1},~~i=0,\ldots,N-1
	\\&&\int \mu_i dv_i = \rho_i,~~i=0,\ldots, N
	\end{eqnarray}
	\end{subequations}
where $KL(\alpha | \beta) = \int \alpha \log\frac{\alpha}{\beta} -\alpha +\beta$ and
	\begin{eqnarray*}
		C_{i,i+1} (x_i,v_i,x_{i+1},v_{j+1}) &=& (12\|x_{i+1} - x_{i} -v_i\|^2 - 12 \langle x_{i+1} - x_{i} -v_i,\, v_{i+1}-v_{i} \rangle 
		\\&&+ 4 \| v_{i+1}-v_{i}\|^2)/(t_{i+1}-t_i).
	\end{eqnarray*}
The optimization variables are joint distributions on the consecutive time points over the phase space. Since the cost is a summation of relative entropies and the constraints are convex, a natural algorithm is that of Bregman projections \cite{BenCarCut15}. 

Define convex constraint sets
	\begin{eqnarray*}
	K_0 &=& \{\int \pi_{01} dx_1dv_1 = \mu_0,~\int \mu_0 dv_0 = \rho_0\},	\\
	K_N &=& \{\int \pi_{N-1,N} dx_{N-1}dv_{N-1} = \mu_N,~\int \mu_N dv_N = \rho_N\}, \mbox{ and}
	\\
	K_i &=& \{\int \pi_{i,i+1} dx_{i+1}dv_{i+1} = \mu_i,~\int \pi_{i-1,i} dx_{i-1}dv_{i-1} = \mu_i,~\int \mu_i dv_i = \rho_i\},
	\end{eqnarray*}
for $i=1,\ldots,N-1$, then the Bregman iterative projection becomes 
	\[
		\pi^n=P_{K_{i_n}}^{KL}(\pi^{n-1}), ~~n=1,2,3,\ldots
	\]
where $i_n$ enumerates $\{0,1,\ldots,N\}$ repeatedly. The projection operator is
	\begin{equation}
		P_{K}^{KL}(\bar \pi):= \argmin_{\pi\in K} KL(\pi \mid \bar\pi).
	\end{equation}
These projections can be derived via Lagrangian method. The projections to $K_0, K_N$ are easier and the projections to $K_1,\ldots,K_{N-1}$ are more involved. Specifically,
	\begin{align*}
		P_{K_0}:\hspace*{.5cm}& \pi_{01} = \frac{\rho_0 \bar\pi_{01}}{\int\bar\pi_{01}dv_0dx_1dv_1}\\
P_{K_N}:\hspace*{.5cm}& \pi_{N-1,N} = \frac{\rho_N \bar\pi_{N-1,N}}{\int \bar\pi_{N-1,N} dv_N dx_{N-1}dv_{N-1}}
	\end{align*}
whereas for $P_{K_i},~i=1,\ldots,N-1$:
	\begin{eqnarray*}
		\pi_{i-1,i} \!&=&\! \frac{\rho_i(\int \bar \pi_{i-1,i} dx_{i-1}dv_{i-1}\int \bar\pi_{i,i+1}dx_{i+1}dv_{i+1})^{1/2}}
		{\int\bar\pi_{i-1,i}dx_{i-1}dv_{i-1}\int(\int \bar \pi_{i-1,i} dx_{i-1}dv_{i-1}\int \bar\pi_{i,i+1}dx_{i+1}dv_{i+1})^{1/2}dv_i}\bar\pi_{i-1,i}
		\\
		\pi_{i,i+1} \!&=&\! \frac{\rho_i(\int \bar \pi_{i-1,i} dx_{i-1}dv_{i-1}\int \bar\pi_{i,i+1}dx_{i+1}dv_{i+1})^{1/2}}
		{\int\bar\pi_{i,i+1}dx_{i+1}dv_{i+1}\int(\int \bar \pi_{i-1,i} dx_{i-1}dv_{i-1}\int \bar\pi_{i,i+1}dx_{i+1}dv_{i+1})^{1/2}dv_i}\bar\pi_{i,i+1}
	\end{eqnarray*}
	
In real implementation, we need to discretize the phase space over a grid.
After discretization, the algorithm only involves matrix multiplication, pointwise-division, multiplication, square root, and therefore can be parallelized easily. The linear convergence rate is guarantee by the property of Bregman projections \cite{BenCarCut15}. Our algorithm should be compared to that developed in \cite{BenGalVia18}. A major difference is that our algorithm doesn't require discretization over the time domain.
	
\section{Conclusion}
We considered a natural extension of the Schr\"odinger bridge problems to multi-marginal partially observable setting. We focused on inertial particles, but more general dynamics can be examined similarly. We discussed the physical meaning, stochastic control formulation and several other aspects of the problems. Just like in the standard SB problem, it has a natural relation to the measure-valued spline theory. An efficient algorithm was also developed, which makes ready for possible applications. We envision that this line of research is going to spark interest in optimal transport theory and application with multiple time points.  


%
%

 \bibliographystyle{splncs04}
 \bibliography{refs}

%
%
%
%
\end{document}